\documentclass[12pt,a4paper]{article}
\usepackage{indentfirst}
\usepackage{amsmath,amssymb,amsfonts,amsthm,graphics}
\usepackage{makeidx}
\usepackage{amsmath,amssymb,amsfonts,amsthm,graphics,graphicx}
\usepackage{makeidx}
\usepackage{color}

\newtheorem{thm}{Theorem}[section]

\newtheorem{lem}{Lemma}[section]
\newtheorem{cor}{Corollary}[section]

\theoremstyle{definition}
\newtheorem{defn}{Definition}[section]

\def\-{\mbox{--}}

\newtheorem{pro}{Proposition}[section]

\def\pf{\noindent {\it Proof.} }

\begin{document}
\title{Proper connection number and \\connected dominating sets\footnote{Supported by NSFC No.11371205,
``973" program No.2013CB834204, and PCSIRT.}}
\author{
\small Xueliang Li, Meiqin Wei, Jun Yue\\
\small Center for Combinatorics and LPMC-TJKLC\\
\small Nankai University, Tianjin 300071, China\\
\small Email: lxl@nankai.edu.cn; weimeiqin8912@163.com; yuejun06@126.com}
\date{}
\maketitle

\begin{abstract}

The proper connection number $pc(G)$ of a connected graph $G$ is
defined as the minimum number of colors needed to color its edges,
so that every pair of distinct vertices of $G$ is connected by at
least one path in $G$ such that no two adjacent edges of the path
are colored the same, and such a path is called a proper path. In
this paper, we show that for every connected graph with diameter 2
and minimum degree at least 2, its proper connection number is 2.
Then, we give an upper bound $\frac{3n}{\delta + 1}-1$ for every
connected graph of order $n$ and minimum degree $\delta$. We also
show that for every connected graph $G$ with minimum degree at least
$2$, the proper connection number $pc(G)$ is upper bounded by
$pc(G[D])+2$, where $D$ is a connected two-way (two-step) dominating
set of $G$. Bounds of the form $pc(G)\leq 4$ or $pc(G)=2$, for many
special graph classes follow as easy corollaries from this result,
which include connected interval graphs, asteroidal triple-free
graphs, circular arc graphs, threshold graphs and chain graphs, all
with minimum degree at least $2$. Furthermore, we get the sharp
upper bound 3 for the proper connection numbers of interval graphs
and circular arc graphs through analyzing their structures.

{\flushleft\bf Keywords}: proper connection number; proper-path
coloring; connected dominating set; diameter

{\flushleft\bf AMS subject classification 2010}: 05C15, 05C40,
05C38, 05C69.

\end{abstract}

\section{Introduction}

All graphs in this paper are finite, connected and simple. An
\emph{edge-coloring} of a graph is a mapping from its edge set to
the set of natural numbers. A path in an edge-colored graph with no
two edges sharing the same color is called a \emph{rainbow path}. An
edge-colored graph $G$ is said to be \emph{rainbow connected} if
every pair of distinct vertices of $G$ is connected by at least one
rainbow path in $G$. Such a coloring is called a \emph{rainbow
coloring} of the graph. The concept of rainbow coloring was
introduced by Chartrand et al. in \cite{CJMZ}. Since then, many
researchers have been studied the problem on the rainbow connection
and got many nice results, see \cite{KY,LS,LS2} for examples. For
more details we refer to a survey paper \cite{LSS} and a book
\cite{LS}.

Inspired by rainbow coloring and proper coloring in graphs, Andrews
et al. \cite{ALLZ} introduced the concept of proper-path coloring.
Let $G$ be an edge-colored graph. A path $P$ in $G$ is called a
\emph{proper path} if no two adjacent edges of $P$ are colored the
same. An edge-coloring $c$ is a \emph{proper-path coloring} of a
connected graph $G$ if every pair of distinct vertices $u,v$ of $G$
are connected by a proper $u$-$v$ path in $G$. If $k$ colors are
used, then $c$ is referred to as a \emph{proper-path $k$-coloring}.
The minimum number of colors needed to produce a proper-path
coloring of $G$ is called the \emph{proper connection number} of
$G$, denoted by $pc(G)$. Form the definition, it is easy to check
that $pc(G)=1$ if and only if $G=K_n$ and $pc(G)=m$ if and only if
$G=K_{1,m}$. For more results, we refer to \cite{ALLZ,BFG}.

A \emph{dominating set} for a graph $G = (V, E)$ is a subset $D$ of
$V$ such that every vertex not in $D$ is adjacent to at least one
member of $D$. The number of vertices in a smallest dominating set
for $G$ is called the \emph{domination number}, denoted by
$\gamma(G)$. The dominating set is every useful to determine some
relationship between a subgraph and its supergraph. There are many
generalized dominating sets, which will be introduced in the
following section and considered in this paper.

We will use two-way dominating sets or a two-way two-step dominating
sets of a graph $G$ to help us find upper bounds of the proper
connection number $pc(G)$. In Section $2$, some definitions and
properties of the proper connection number of a graph are given. In
Section $3$, we give the bound $pc(G)\leq pc(G[D])+2$, where $D$ is
a connected two-way two-step dominating set of $G$. And we get the
following two results as its corollaries: One is that the proper
connection number of a chain graph with minimum degree at least 2 is
$2$; the other is that for every connected graph of order $n$ and
minimum degree $\delta$, its proper connection number is upper
bounded by $\frac{3n}{\delta+1}-1$. In addition, we also get that a
graph with diameter 2 and minimum degree at least 2 has proper
connection number 2. In Section $4$, we turn to using connected
two-way dominating sets $D$ of $G$. The inequality $pc(G)\leq
pc(G[D])+2$ and upper bounds for interval graphs, asteroidal
triple-free graphs, circular arc graphs and threshold graphs are
obtained. Furthermore, we get the sharp upper bound 3 for proper
connection numbers of interval graphs and circular arc graphs
through analyzing their structures.

\section{Preliminaries}

In this section, we introduce some definitions and present several
useful facts about the path connection numbers of graphs. We begin
with some basic conceptions.
\begin{defn}
Let $G$ be a connected graph. The \emph{distance between two
vertices} $u$ and $v$ in $G$, denoted by $d(u,v)$, is the length of
a shortest path between them in $G$. The \emph{eccentricity} of a
vertex $v$ is $ecc(v):=max_{x\in V(G)}d(v, x)$. The \emph{diameter}
of $G$ is $diam(G):=max_{x\in V(G)}ecc(x)$. The \emph{radius} of $G$
is $rad(G):=min_{x\in V(G)}ecc(x)$. The \emph{distance between a
vertex $v$ and a set $S\subseteq V(G)$} is $d(v,S):=min_{x\in S}d(v,
x)$. The \emph{$k$-step neighborhood} of a set $S\subseteq V(G)$ is
$N^k(S):=\{x\in V(G)|d(x,S)=k\}, k\in \{0, 1, 2, \ldots\}$. The
\emph{degree} of a vertex $v$ is $deg(v):=|N^1({v})|$. The
\emph{minimum degree} of G is $\delta(G):=min_{x\in V(G)}deg(x)$. A
vertex is called \emph{pendant} if its degree is $1$ and
\emph{isolated} if its degree is $0$. We may use $N^k(v)$ in place
of $N^k(\{v\})$.
\end{defn}
\begin{defn}
Given a graph $G$, a set $D\subseteq V(G)$ is called a
\emph{$k$-step dominating set} of $G$, if every vertex in $G$ is at
a distance at most $k$ from $D$. Further, if $D$ induces a connected
subgraph of $G$, it is called a {\it connected $k$-step dominating
set} of $G$.
\end{defn}
\begin{defn}
A dominating set $D$ in a graph $G$ is called a \emph{two-way
dominating set}, if every pendant vertex of G is included in $D$. In
addition, if $G[D]$ is connected, we call $D$ a \emph{connected
two-way dominating set}.
\end{defn}
\begin{defn}
A two-step dominating set $D$ of vertices in a graph $G$ is called a
\emph{two-way two-step dominating set} if

(i) every pendant vertex of $G$ is included in $D$ and

(ii) every vertex in $N^2(D)$ has at least two neighbors in
$N^1(D)$.

Further, if $G[D]$ is connected, $D$ is called a {|it connected
two-way two-step dominating set} of $G$.
\end{defn}
\begin{defn}
Let $G$ be a graph and $s$ a positive integer. Define $sG$ as
disjoint union of $s$ copies of the graph $G$, i.e.,
$sG=\underbrace{G\cup G\cup\cdots\cup G}_s.$
\end{defn}
\begin{defn}
A \emph{Hamiltonian path} in a graph $G$ is a path containing every
vertex of $G$. And a graph having a Hamiltonian path is called a
\emph{traceable graph}.
\end{defn}

We state some known simple results on the proper-path coloring and
two-way two-step dominating set which will be useful in the sequel..
\begin{lem}\label{lem2.1}
If $P$ is a path, then $pc(P)=2$.
\end{lem}
\begin{lem}\cite{ALLZ}\label{lem2.3}
If $G$ is a traceable graph that is not complete, then $pc(G) = 2$.
\end{lem}
\begin{pro}\cite{ALLZ}\label{prop1}
If $T$ be a nontrivial tree, then $pc(T) = \chi' = \Delta$.
\end{pro}

In the same paper \cite{ALLZ}, there is a lemma which will be useful
in the following proof.
\begin{lem}\cite{ALLZ}\label{lem2.2}
If $G$ is a nontrivial connected graph and $H$ is a connected
spanning subgraph of $G$, then $pc(G)\leq pc(H)$. In particular,
$pc(G)\leq pc(T)$ for every spanning tree $T$ of $G$.
\end{lem}

\section{Proper connection number and connected two-way two-step dominating set}

In this section, we will give an upper bound of proper connection
number of a graph $G$ by using the connected two-way two-step
dominating sets.

Let $D$ be a connected two-way two-step dominating set of a graph
$G$. This implies that every vertex $v\in V(G)\setminus D$ has two
edge-disjoined paths connecting to $D$. Our idea is to color $G[D]$
first and then to color all the other edges with a constant number
of colors, ensuring a proper-path coloring of a graph $G$. Now we
give our main theorems.
\begin{thm}\label{thm3.1}
If $D$ is a connected two-way two-step dominating set of a graph $G$, then
$$
pc(G)\leq pc(G[D])+2.
$$
\end{thm}

\pf Let $H$ be a spanning subgraph of a graph $G$. Then $pc(G)\leq
pc(H)$ by Lemma \ref{lem2.2}. In the following, we will give a
proper-path coloring of $H$ with $pc(G[D])+2$ colors, and then prove
the theorem.

Let $c_D$ be a proper-path coloring of $G[D]$ using colors
$\{3,4,\cdots, k:=pc(G[D]), k+1,k+2\}$. For $x\in N^1(D)$, a
neighbor of $x$ in $D$ is called a foot of $x$. Define the set of
foots of $x$ as $F(x)=\{u: u\ \text{is a foot of}\ x\}$. And define
the set of the neighbors of a vertex $v \in N^2(D)$ in $N^1(D)$ to
be $F^1(v) = \{u : u~ \text{is the neighbor of}~ v~ \text{in}
~N^1(D)\}$.

\textbf{Case 1.} For each vertex $v\in N^2(D)$, its neighbors in
$N^1(D)$ has at least one common foot. That is to say, the set
$N^2(D)=\{v_1,v_2,\cdots,v_t\}$ and the neighborhood of $v_i$ in
$N^1(D)$ is $F^1(v_i)=\{u_{i,1},u_{i,2},\cdots,u_{i,i_l}\}$, where
$|F^1(v_i)|\geq 2,\ (i=1,2,\cdots,t)$. Then
$\bigcap\limits_{a=1}^{i_l}F(u_i,a)\neq \emptyset\
(i=1,2,\cdots,t)$.

\begin{figure}[h,t,b,p]
\begin{center}
\includegraphics[bb = 112 515 596 715,scale = 0.8]{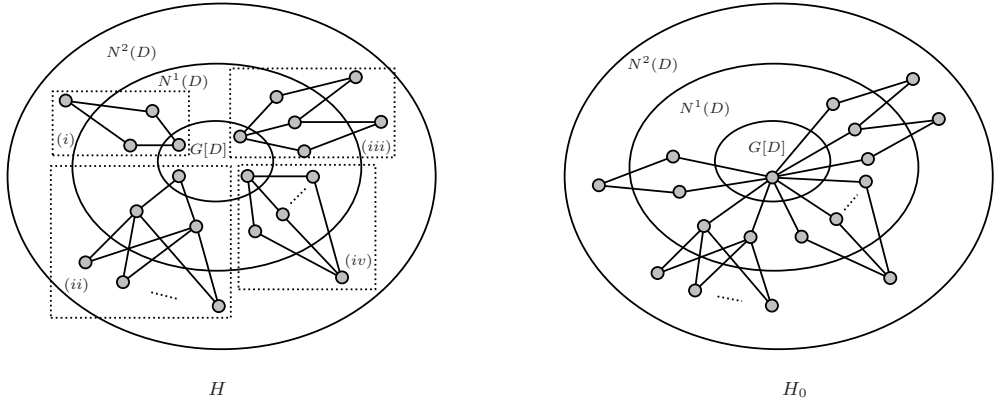}
\caption{Spanning subgraphs of $G$} \label{fig1}
\end{center}
\end{figure}

In this case, $p(i)\cup q(ii)\cup r(iv)\cup s(iv)\cup G[D]$ (see
Figure \ref{fig1}, where the $(i),~(ii),~(iv),~(vi)$ are the
subgraphs and $p,~q,~r,~s$ are the numbers of the corresponding
subgraphs in $G$) is a spanning subgraph of $G$, in which we do not
exclude the case that the foots of some vertices are in common.
Since each pair of vertices $x,y\in D$ has a proper $x$-$y$ path in
$G[D]$ under the coloring $c_D$, it suffices to show that $p(i)\cup
q(ii)\cup r(iv)\cup s(iv)\cup G[D]$, in which all the vertices in
$N^1(D)$ have one common root, has a proper-path coloring using
$k+2$ distinct colors.

\begin{figure}[h,t,b,p]
\begin{center}
\includegraphics[bb= 121 516 593 715,scale = 0.8]{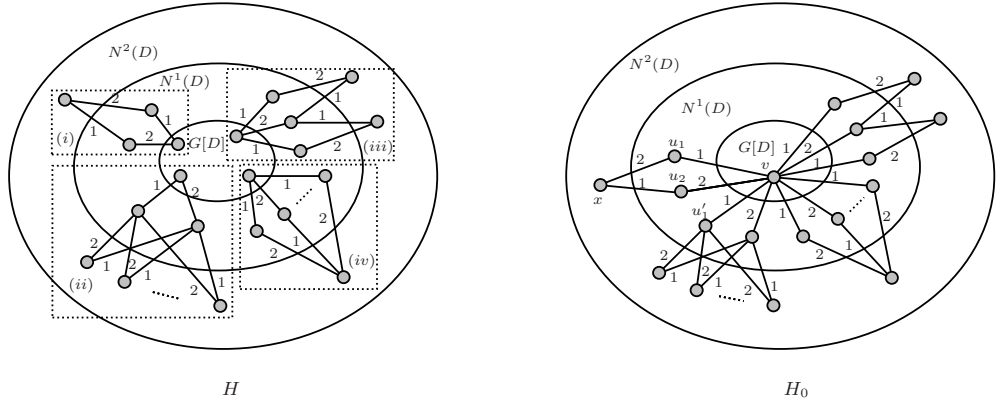}
\caption{The proper-path coloring for the spanning subgraphs of $G$} \label{fig2}
\end{center}
\end{figure}

Give an edge-coloring $c$ using colors $\{1,2,\cdots,k,k+1,k+2\}$
for the above spanning subgraphs $H,H_0$ of $G$ as follows: for the
edges in $G[D]$, we use the proper-path coloring $c_D$; and for the
edges in $(i),(ii),(iii)$ and $(iv)$, color them as depicted in
Figure \ref{fig2}. Then for any two vertices $u_i,u'_i\in N^1(D)$,
we can find a proper $u_i$-$u'_i$ path as follows: if
$c(u_iv)=c(u'_iv)$, then $u_ixu_jvu'_i$ is a proper $u_i$-$u'_i$
path; if $c(u_iv)\neq c(u'_iv)$, then $u_ivu'_i$ is a proper
$u_i$-$u'_i$ path. For every pair of vertices $u,v\in N^2(D)$ or
$u\in N^1(D),v\in N^2(D)$ or $u\in N^1(D)\cup N^2(D),v\in D$, there
exist a proper $u$-$v$ path under the coloring $c$ as well. This
implies that $c$ is a proper-path coloring of the graph $H_0$ and it
follows that $pc(G)\leq pc(H_0)\leq pc(G[D])+2$.

\textbf{Case 2.} There exists one vertex $x\in N^2(D)$ whose
neighbors in $N^1(D)$ has no common roots. Note that such vertices
are not necessarily unique and we can similarly prove the same
result as it for $x$ in this theorem.

We give a proper-path coloring $c$ using colors
$\{1,2,\cdots,k:=pc(G[D]),k+1,k+2\}$ for a spanning subgraphs $H$ of
$G$ as well. Similarly, for the edges in $G[D]$, we still use the
proper-path coloring $c_D$. By the definition of connected two-way
two-step dominating sets, $x$ has at least two distinct neighbors in
$N^1(D)$ and two edge-disjoined paths connecting to $D$. This
implies that there exist two vertex-disjoint paths, denoted by
$P_1=xu_iv_i,P_2=xu_jv_j$, where $u_i,u_j\in N^1(D)$ and $v_i,v_j\in
D$. We color the edges $xu_i$ with color $1$ or color $2$ such that
$\{1,2\}\subseteq\{c(xu_i):u_i\in N^1(D)\}$ holds for every vertex
$x\in N^2(D)$. And set $c(u_iv_i)\in \{1,2\}\setminus c(xu_i)$. Then
for any two vertices $u_i,u'_i\in N^1(D)$, we can find a proper
$u_i$-$u'_i$ path as follows: if $v_i\neq v'_i$, then
$u_iv_iP_{ii'}v'_iu'_i$ is a proper $u_i$-$u'_i$ path, in which
$P_{ii'}$ is a proper $v_i$-$v'_i$ path in $G[D]$; if $v_i=v'_i$ and
$c(u_iv_i)=c(u'_iv'_i)$, then $u_ixu_jv'_iu'_i$ is a proper
$u_i$-$u'_i$ path, where $u_j$ is a neighbor of $x$ in $N^1(D)$ such
that $c(u_ix)\neq c(xu_j)$; if $v_i=v'_i$ and $c(u_iv_i)\neq
c(u'_iv'_i)$, then $u_iv_iu'_i$ is a proper $u_i$-$u'_i$ path.
Similarly, one can check that for every pair of vertices $u,v\in
N^2(D)$ or $u\in N^1(D),v\in N^2(D)$ or $u\in N^1(D)\cup N^2(D),v\in
D$, there exists a proper $u$-$v$ path under the coloring $c$. It
means that $c$ is a proper-path coloring of a spanning subgraph $H$
of $G$, and then $pc(G)\leq pc(H)\leq pc(G[D])+2$.
\begin{figure}[h,t,b,p]
\begin{center}
\includegraphics[bb=160 530 433 704,scale = 0.8]{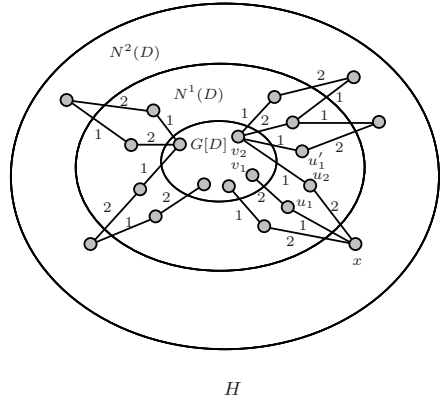}
\caption{An example for the proper-path coloring of the spanning subgraph} \label{fig3}
\end{center}
\end{figure}
\vskip -5mm
\qed

\begin{figure}[h,t,b,p]
\begin{center}
\includegraphics[bb=215 600 368 716,scale = 0.8]{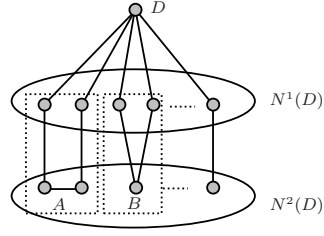}
\caption{Two structure of the $N^2(D)$ adjacent to $N^1(D)$} \label{fig4}
\end{center}
\end{figure}

Let $G$ be a graph with diameter 2 and minimum degree at least 2.
Then there exists one vertex in $G$ which forms a two-way two-step
dominating set. And the vertices in $N^2(D)$ must be contained in
the following two structures $A$ and $B$ (see Figure \ref{fig4}),
since the diameter of $G$ is 2 and minimum degree is at least 2.

Thus, by the proof of the Theorem \ref{thm3.1}, we can easily get
the following corollary. \vskip -5mm
\begin{cor}
Let $G$ be a graph with diameter 2 and minimum degree at least 2.
Then $pc(G) = 2$.
\end{cor}

\begin{figure}[h,t,b,p]
\begin{center}
\includegraphics[bb=152 603 394 733,scale = 0.8]{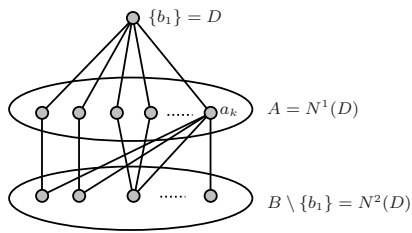}
\caption{Graph for Corollary 3.2} \label{fig5}
\end{center}
\end{figure}

A bipartite graph $G(A,B)$ is called a \emph{chain graph}, if the
vertices of $A$ can be ordered as $A=(a_1,a_2,\ldots,a_k)$ such that
$N(a_1)\subseteq N(a_2)\subseteq\cdots\subseteq N(a_k)$. Applying
Theorem \ref{thm3.1}, we can give the proper connection number of a
connected chain graph as follows.
\begin{cor}
If $G$ is a connected chain graph with minimum degree at least 2,
$pc(G)=2$.
\end{cor}
\pf Let $G=G(A,B)$ be a connected chain graph, where
$A=(a_1,a_2,\ldots,a_k),\ B=(b_1,b_2,\ldots,b_s)$ such that $b_1\in
N(a_1)\subseteq N(a_2)\subseteq\cdots\subseteq N(a_k)$. Obviously,
$N(a_k)=B$ since $G=G(A,B)$ is connected. It is easy to verify that
$D=\{b_1\}$ is a connected two-way two-step dominating set of $G$
and $N^1(D)=A,\ N^2(D)=B\setminus \{b_1\}$ (see Figure \ref{fig5}).
Applying the result in Theorem \ref{thm3.1}, we obtain that
$pc(G)\leq 2$. On the other hand, $pc(G)=1$ if and only if $G=K_n$
and then $pc(G)\geq 2$. Therefore, $pc(G)=2$.

\qed

In \cite{CDR}, there is a lemma giving the size of a connected
two-way two-step dominating set, which is stated as follows.
\begin{lem}\cite{CDR}\label{lem2.4}
Every connected graph $G$ of order $n\geq 4$ and minimum degree
$\delta$ has a connected two-way two-step dominating set $D$ of size
at most $\frac{3n}{\delta+1}-2$.
\end{lem}

Then, we get the following corollary.
\begin{cor}
Let $G$ be a connected graph of order $n\geq 4$ and minimum degree $\delta$.
Then we have
$$
pc(G)\leq \frac{3n}{\delta+1}-1.
$$
\end{cor}
\pf By Lemma \ref{lem2.4}, the graph $G$ has a connected two-way
two-step dominating set $D$ such that $|D|\leq
\frac{3n}{\delta+1}-2$. Since $G[D]$ is connected, $pc(G[D])\leq
|D|-1\leq \frac{3n}{\delta+1}-3$. Together with the result in
Theorem \ref{thm3.1}, we obtain that
$$
pc(G)\leq pc(G[D])+2\leq \frac{3n}{\delta+1}-1.
$$
\qed

{\bf Remark 1.} This upper bound of proper connection number is not
sharp. Further effort is needed to find a sharp upper bound.

{\bf Remark 2.} If the minimum degree of a graph is at least
$\frac{n}{2}$, then the graph is Hamiltonian, and then $pc(G) =2$.
But the corollary shows that if there exists some $k$ such that
$\delta=kn$, then the proper connection number can be upper bounded
by $\frac{3}{k}-1$, where $k \leq \frac{1}{2}$.

{\bf Remark 3.} Since $G[D]$ is a connected subgraph of a graph $G$,
by Proposition \ref{prop1} we know that $pc(G) \leq \chi'(T) + 2$,
where $T$ is a spanning tree of $G[D]$.

\section{Proper connection number and connected two-way dominating set}

The definition of a two-way dominating set ($D$) implies that every
vertex in $V(G)\setminus D$ has at least two edge-disjoint paths
connecting to $D$. Similar with the idea in Section $3$, we also
obtain the following upper bound for proper connection number.
\begin{thm}\label{thm4.1}
If $D$ is a connected two-way dominating set of a graph $G$, then
$$
pc(G)\leq pc(G[D])+2.
$$
\end{thm}
\pf We will give a proper-path coloring of a spanning subgraph $H$
of the graph $G$ with $pc(G[D])+2$ colors, which implies this
theorem. Let $c_D$ be a proper-path coloring of $G[D]$ using colors
$\{3,4,\cdots, k:=pc(G[D]), k+1,k+2\}$.

\begin{figure}[h,t,b,p]
\begin{center}
\includegraphics[bb=105 578 483 727,scale = 0.8]{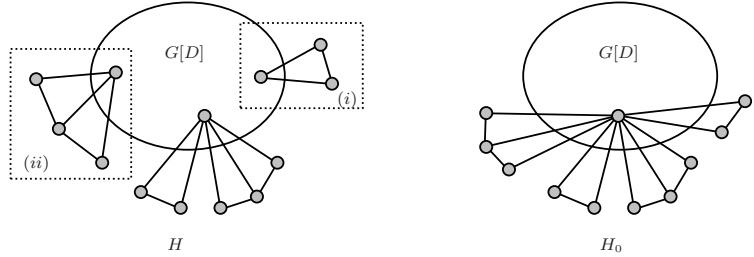}
\caption{Spanning subgraphs of $G$} \label{fig6}
\end{center}
\end{figure}

For any $x\in G\setminus D$, we call a neighbor of $x$ in $D$ a foot
of $x$. Define the set of the foots of $x$ as $F(x)=\{u: u\ \text{is
a foot of}\ x\}$. We focus on the case that $|F(x)|=1$ for every
vertex $x\in G\setminus D$. Since $D$ is a connected two-way
dominating set, every pendant vertex of $G$ is included in $D$.
Additionally, each pair of vertices $x,y\in D$ has a proper $x$-$y$
path in $G[D]$ under the coloring $c_D$ and two colors are enough to
ensure that a path is proper. Consequently, $p(i)\cup q(ii)\cup
G[D]$ (see Figure \ref{fig6}) is a spanning subgraph of $G$, where
we allow that the foots of some vertices are in common. It suffices
to show that $p(i)\cup q(ii)\cup G[D]$, in which all the vertices in
$G\setminus D$ have a common root, has a proper-path coloring using
$k+2$ distinct colors.

\begin{figure}[h,t,b,p]
\begin{center}
\includegraphics[bb=90 579 490 729,scale = 0.8]{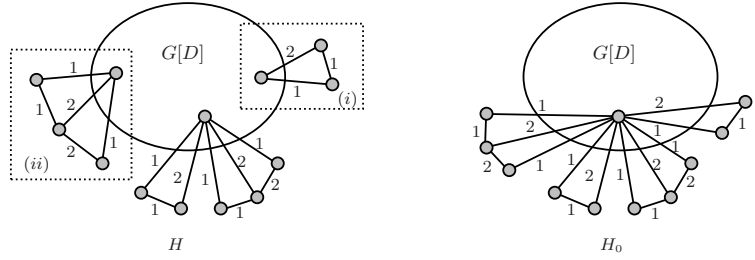}
\caption{The proper-path coloring for the spanning subgraphs of $G$} \label{fig7}
\end{center}
\end{figure}

Now we give an edge-coloring $c$ using colors
$\{1,2,\cdots,k,k+1,k+2\}$ for the above spanning subgraphs $H,H_0$
of $G$ as follows: for the edges in $G[D]$, we use the proper-path
coloring $c_D$; and for the edges in $(i),(ii)$, color them as
depicted in Figure \ref{fig7}. Then for any two vertices $u_i,u_j\in
G\setminus D$, a proper $u_i$-$u_j$ path can be found in $H_0$ under
the coloring $c$ as follows: if $u_iu_j\in H_0$, then $u_iu_j$ is a
proper $u_i$-$u_j$ path; if $u_iu_j\notin H_0$, $u_i\in (i)$ and
$u_j\in (i)$, then $u_ivu_j$ or $u_ivu_ku_j\ (u_ku_j\in H_0)$ is a
proper $u_i$-$u_j$ path; if $u_iu_j\notin H_0$, $u_i\in (ii)$ and
$u_j\in (ii)$, then $u_ivu_j$ or $u_ivu_ku_j\ (u_ku_j\in H_0)$ or
$u_iu_kvu_j\ (u_iu_k\in H_0)$ is a proper $u_i$-$u_j$ path; if
$u_iu_j\notin H_0$, $u_i\in (i)$ and $u_j\in (ii)$, then $u_ivu_j$
or $u_iu_kvu_j\ (u_iu_k\in H_0)$ is a proper $u_i$-$u_j$ path. One
can find a proper $u$-$v$ path for every pair of vertices $u\in D$
and $v\in G\setminus D$ in similarly way. This implies that $c$ is a
proper-path coloring of the graph $H_0$, and then $pc(G)\leq
pc(H_0)\leq pc(G[D])+2$.

Next we consider that there exist some vertices $x_i\in G\setminus
D$ such that $|F(x_i)|\geq 2$. Let $u_{i_1},u_{i_2}\in F(x_i)$ for
every such vertex $x_i$. On the basis of the coloring in the above
case, color $u_{i_1}x_i,u_{i_2}x_i$ with color $1$ if they have not
been colored. This provides a proper path for $x_i$ and every other
vertices in $G$.

Summarizing the above analysis, this theorem holds. \qed

As consequences of Theorem\ref{thm4.1}, the upper bounds for proper
connection numbers of interval graphs, asteroidal triple-free
graphs, circular arc graphs, threshold graphs, and chain graphs are
followed. Before presenting the upper bounds, we state the
definitions of all these graphs first.
\begin{defn}
An \emph{intersection graph} of a family of sets F is a graph whose
vertices can be mapped to sets in $F$ such that there is an edge
between two vertices in the graph if and only if the corresponding
two sets in $F$ have a nonempty intersection. An \emph{interval
graph} is an intersection graph of intervals on the real line. A
\emph{circular arc graph} is an intersection graph of arcs on a
circle.
\end{defn}
\begin{defn}
An independent triple of vertices $x,y,z$ in a graph $G$ is an
\emph{asteroidal triple (AT)}, if between every pair of vertices in
the triple, there is a path that does not contain any neighbor of
the third. A graph without asteroidal triples is called an
asteroidal \emph{triple-free (AT-free) graph}.
\end{defn}
\begin{defn}
A graph $G$ is a \emph{threshold graph}, if there exists a weight
function $w: V(G)\rightarrow R$ and a real constant $t$ such that
two vertices $u, v\in V(G)$ are adjacent if and only if
$w(u)+w(v)\geq t$.
\end{defn}

For a graph $G$ with $\delta(G)\geq 2$, every (connected) dominating
set of $G$ is a (connected) two-way dominating set. Next, we will
give some upper bounds for the proper connection numbers of the
above classes of graphs.
\begin{cor}\label{cor4.2}
Let $G$ be a connected non-complete graph with $\delta(G)\geq 2$. Then

~~(i) if $G$ is an interval graph, $pc(G)\leq 4$,

~(ii) if $G$ is AT-free, $pc(G)\leq 4$,

(iii) if $G$ is a circular arc graph, $pc(G)\leq 4$,

(iv) if $G$ is a threshold graph, $pc(G)=2$,
\end{cor}

{\bf Remark.} There are four well-known results on the dominating
sets of the graphs in Corollary \ref{cor4.2}, which are stated as
follows: (i) every interval graph $G$ which is not isomorphic to a
complete graph has a dominating path of length at most $diam(G)-2$,
(ii) every AT-free graph $G$ has a dominating path of length at most
$diam(G)$, (iii) every circular arc graph $G$, which is not an
interval graph, has a dominating cycle of diameter at most
$diam(G)$, (iv) a maximum weight vertex in a connected threshold
graph $G$ is a dominating vertex. Together with Theorem \ref{thm4.1}
and Lemma \ref{lem2.1}, we obtain the upper bounds in Corollary
\ref{cor4.2}.

Furthermore, we can get a better and sharp bounds on their proper
connection numbers by analyzing the structure of the connected
interval graph and circular arc graph, which can be stated as
follows.
\begin{thm}\label{thm4.2}
Let $G$ be a connected interval or circular arc graph with
$\delta(G)\geq 2$. Then the proper connection number $pc(G)\leq 3$,
and this bound is sharp.
\end{thm}
\pf The proofs for connected interval graph and circular arc graph
are similar, since the circular arc graph is a generalization. So we
only give the details for $G$ being a connected interval graph with
minimum degree at least 2. We will give a proper-path coloring of
$G$ using colors $1,2,3$ and as well some graphs arriving at this
bound.
\begin{figure}[h,t,b,p]
\begin{center}
\includegraphics[bb=227 518 556 709,scale = 0.8]{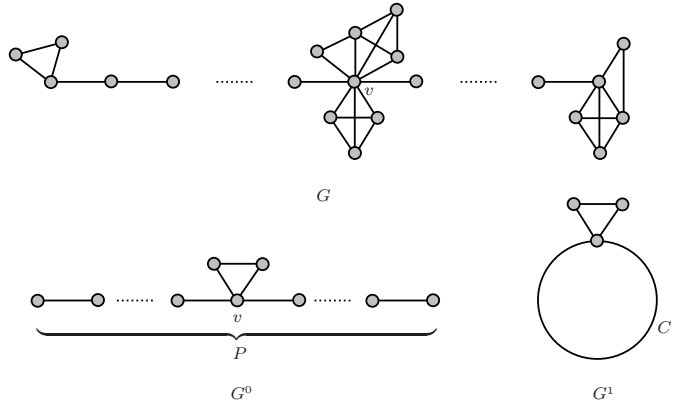}
\caption{Graphs for Theorem 4.2} \label{fig8}
\end{center}
\end{figure}

Let $P$ be a dominating path of length of $diam(G)-2$ in $G$. We
color the edges of $P$ using colors 1 and 2, alternately.
Additionally, by the definition of an interval graph, for every
vertex $v\in P$ the subgraph of $G[G\setminus P\cup\{v\}]$
containing $v$ is a maximal clique having at least one common vertex
$v$ (see Figure \ref{fig8}). For convenience, we call the vertex $v$
the root of those maximal cliques. For a fixed vertex $v\in P$, we
give a coloring for the edges in those maximal cliques containing
$v$ and for any other vertices the same. Let $Q_1,Q_2,\cdots,Q_t$
denote the maximal cliques whose common root is $v$. Color
$u_1v,u_2v,\cdots,u_tv$ with color 3, where $u_i\in Q_i$ and
$u_1,u_2,\cdots,u_t$ are not necessarily different. And color all
the other edges in $Q_1,Q_2,\cdots,Q_t$ with color 1. One can check
that such a coloring is a proper-path coloring of $G$.

Actually, there are many connected interval graphs (or circular arc
graphs) with proper connection number $3$, which implies that this
bound is sharp. Here we give an infinite class of such graphs. As
depicted in Figure \ref{fig8}, $G_0 \ ( \text{or}~G_1)$, for any
vertex $v\in P(\text{or}~C)$, the subgraph of $G[G\setminus
P(\text{or}~ C)\cup\{v\}]$ containing $v$ is a triangle. And for any
other vertices in $P(\text{or}~C)$, those cliques are arbitrary. It
is easy to verify that two colors are not enough to make the
coloring proper for $G$. \qed

\end{document}